\documentclass[10pt]{article}
\usepackage{theorem,amsfonts,amsmath,latexsym,mathrsfs,amssymb,upref}

\newcommand{\D}{{\mathbf D}}
\newcommand{\N}{{\mathbf N}}

\newcommand{\R}{{\mathbf R}}

\newcommand{\Z}{{\mathbf Z}}
\newcommand{\C}{{\mathbf C}}
\newcommand{\rmP}{{\rm P}}
\newcommand{\rmR}{{\rm R}}

\newtheorem{thm}{Theorem}
\newtheorem{cor}{Corollary}

\allowdisplaybreaks

\begin{document}

{\noindent Generalizations and specializations of 
generating functions for Jacobi, Gegenbauer, Chebyshev and Legendre polynomials
with definite integrals\\[-0.2cm]}

%

\noindent Howard S.\ Cohl$^\ast$ and Connor MacKenzie$^\dag$\\[0.2cm]
$^\ast$~Applied and Computational Mathematics Division,
National Institute of Standards and Technology,
Gaithersburg, Maryland, U.S.A.\\[-0.2cm]


\noindent $^\dag$~Department of Mathematics, Westminster College, 319 South Market Street, 
New Wilmington, Pennsylvania, U.S.A.


\begin{abstract}
\noindent 
In this paper we generalize and specialize generating functions for 
classical orthogonal
polynomials, namely Jacobi, Gegenbauer, Chebyshev and Legendre polynomials.
We derive a generalization of the generating function 
for Gegenbauer polynomials through extension a two element sequence 
of generating functions for Jacobi polynomials.
Specializations of generating functions are accomplished through the
re-expression of Gauss hypergeometric functions in terms of less
general functions.  Definite integrals which correspond to the presented 
orthogonal polynomial series expansions are also given.
\end{abstract}






\newpage
\section{Introduction}
\label{Introduction}

This paper concerns itself with analysis of generating functions
for Jacobi, Gegenbauer, Chebyshev and Legendre polynomials.  The analysis
involves generalization and specialization by re-expression of 
Gauss hypergeometric generating functions for these orthogonal
polynomials.
The generalizations that we present here are for two of the most 
important generating functions for Jacobi polynomials, namely
\cite[(4.3.1-2)]{Ismail}.\footnote{The interesting question 
of which orthogonal polynomial generating functions are important is addressed in 
\cite[Section 4.3]{Ismail} (see also R.~A.~Askey's MathSciNet review of 
Srivastava \& Manocha (1984) \cite{SriManocha}).} In fact, these are
the first two generating functions which appear in Section 4.3 of \cite{Ismail}.
As we will show, these two generating functions, traditionally expressed in 
terms of Gauss hypergeometric functions, can be re-expressed in terms 
of associated Legendre functions (and also in terms of Ferrers functions,
associated Legendre functions on the real segment $(-1,1)$).  
Our Jacobi polynomial generating function 
generalizations, Theorem \ref{firstbigthm}, Corollary \ref{firstbigcor} and 
Corollary \ref{Szegocorollary}, generalize the generating function
for Gegenbauer polynomials.
The presented proofs of these generalizations rely upon 
the series re-arrangment technique.
The motivation for the proofs of our 
generalizations was purely intuitive.  Examination of the two
Jacobi polynomial generating functions which we generalize,\footnote{As well 
as their companion identities, see Section 
\ref{Jacobifunctiongeneralization}.} indicate 
that these generating functions represent two elements of an infinite
sequence of eigenfunction expansions.  The resulting formal proofs of 
our generalizations are then simple consequences.

Our generalized expansions and hypergeometric orthogonal polynomial
generating functions are given in terms of Gauss hypergeometric functions.
The Gauss hypergeometric Jacobi polynomial generating functions 
which we generalize, as well as their eigenfunction expansion generalizations, 
are all re-expressible in terms of associated Legendre functions.  
Associated Legendre functions \cite[Chapter 14]{NIST} are given 
in terms of Gauss hypergeometric functions which satisfy a quadratic 
transformation of variable.  These have an
abundance of applications in Physics, Engineering and Applied Mathematics 
for solving partial differential equations in a variety of contexts.
Recently, efficient numerical evaluation of these functions has been 
investigated in \cite{Schneideretal}.
Associated Legendre functions are more elementary than Gauss hypergeometric functions
because Gauss hypergeometric functions have three free parameters, whereas associated 
Legendre functions have only two. One can make this argument of elementarity with 
all of the functions which can be expressed in terms of Gauss hypergeometric functions,
namely (inverse) trigonometric, (inverse) hyperbolic, exponential, logarithmic, 
Jacobi, Gegenbauer, and Chebyshev polynomials, and complete elliptic integrals.
We summarize how associated Legendre functions, Gegenbauer, Chebyshev 
and Legendre polynomials and complete elliptic integrals of the 
first kind are interrelated.  
In a one-step process, we obtain definite integrals from 
our orthogonal polynomial expansions and generating functions. 

To the best of our knowledge our generalizations, re-expressions
of Gauss hypergeometric generating functions for orthogonal polynomials 
and definite integrals are new and have not previously 
appeared in the literature.  Furthermore,
the generating functions presented in this paper are some of the most
important generating functions for these hypergeometric orthogonal 
polynomials and any
specializations and generalizations will be similarly important.

This paper is organized as follows.  In Sections
\ref{Jacobifunctiongeneralization},
\ref{Gegenbauerpolynomials},
\ref{Chebyshevpolynomialsofthesecondkind},
\ref{Legendrepolynomialsofthefirstkind},
\ref{Chebyshevpolynomialsofthefirstkind}, we present
generalized and simplified expansions for Jacobi, Gegenbauer,
Chebyshev of the second kind, Legendre,
and Chebyshev of the first kind polynomials respectively.
In \ref{DefiniteintegralsforJacobipolynomials} we present 
definite integrals which correspond to the derived hypergeometric
orthgonal polynomial expansions.
Unless stated otherwise the domains of convergence given
in this paper are those of the original generating function and/or
its corresponding definite integral.

\medskip

\noindent Throughout this paper we rely on the following definitions.  
Let $a_1,a_2,a_3,\ldots\in\C$. If $i,j\in\Z$ and $j<i,$ then
$\sum_{n=i}^{j}a_n=0$ and $\prod_{n=i}^ja_n=1$.
The set of natural numbers is given by $\N:=\{1,2,3,\ldots\}$,
the set $\N_0:=\{0,1,2,\ldots\}=\N\cup\{0\}$, and 
$\Z:=\{0,\pm 1,\pm 2,\ldots\}.$
Let $\D:=\{z\in\C:|z|<1\}$ be the open unit disk.

\section{Expansions over Jacobi polynomials}
\label{Jacobifunctiongeneralization}

The Jacobi polynomials $P_n^{(\alpha,\beta)}:\C\to\C$ can be defined 
in terms of the terminating Gauss hypergeometric series as follows
(\cite[(18.5.7)]{NIST})
\[
P_n^{(\alpha,\beta)}(z):=
\frac{(\alpha+1)_n}{n!}\,
{}_2F_1\left(
\begin{array}{c}
-n,n+\alpha+\beta+1\\[0.1cm]
\alpha+1
\end{array};
\frac{1-z}{2}
\right),
\]
for $n\in\N_0$, and $\alpha,\beta>-1$ such that if $\alpha,\beta\in(-1,0)$ then 
$\alpha+\beta+1\ne 0$.
The Gauss hypergeometric function
${_2}F_1:\C^2\times(\C\setminus\N_0)\times\D\to\C$ 
(see Chapter 15 in 
\cite{NIST}) is defined as
\[
{_2}F_1\left(
\begin{array}{c}
a,b\\[0.2cm]
c
\end{array};z\right)
:=\sum_{n=0}^\infty \frac{(a)_n(b)_n}{(c)_n} \frac{z^n}{n!},
\]
\noindent where the Pochhammer symbol (rising factorial) $(\cdot)_n:\C\to\C$ \cite[(5.2.4)]{NIST} is defined by
\[
(z)_n:=\prod_{i=1}^n(z+i-1),\nonumber
\]
where $n\in\N_0$. 
Note that the Gauss hypergeometric function can be analytically continued 
through for instance the Euler's integral representation
for $z\in\C\setminus(1,\infty)$ (see for instance 
\cite[Theorem 2.2.1]{AAR}).

Consider the generating function for Gegenbauer polynomials
(see \S\ref{Gegenbauerpolynomials} for their definition)
given by 
\cite[(18.12.4)]{NIST}, namely
\begin{equation}
\frac{1}{(1+\rho^2-2\rho x)^\nu}=\sum_{n=0}^\infty\rho^n C_n^\nu(x).
\label{generatingfunctionforGeg}
\end{equation}
We attempt to generalize this expansion using
the representation of Gegenbauer polynomials in terms of Jacobi polynomials
given by 
\cite[(18.7.1)]{NIST}, namely
\begin{equation}
C_n^\nu(x)=\frac{(2\nu)_n}{\left(\nu+\frac12\right)_n}P_n^{(\nu-1/2,\nu-1/2)}(x).
\label{JactoGeg}
\end{equation}
By making the replacement $\nu-1/2$ to $\alpha$ and $\beta$ in
(\ref{generatingfunctionforGeg}) using (\ref{JactoGeg}), we see that
there are two possibilities for generalizing the generating 
function for Gegenbauer polynomials to a generating function for Jacobi polynomials.
These two possibilities are given below, namely 
(\ref{genGeg1}),
(\ref{genGeg1companion}).
The first possibility is given for $\rho\in\D\setminus(-1,0]$ by
\cite[(18.12.3)]{NIST}
\begin{eqnarray}
&&\hspace{-1.2cm}\frac{1}{(1+\rho)^{\alpha+\beta+1}}
\,{}_2F_1\left(
\begin{array}{c}
\frac{\alpha+\beta+1}{2},\frac{\alpha+\beta+2}{2}
\\[0.1cm]
\beta+1
\end{array};
\frac{2\rho(1+x)}{(1+\rho)^2}
\right)\nonumber\\[0.2cm]
&&\hspace{-0.2cm}=\left(\frac{2}{\rho(1+x)}\right)^{\beta/2}
\frac{\Gamma(\beta+1)}
{\rmR^{\alpha+1}}
P_\alpha^{-\beta}
(\zeta_+)
=\sum_{n=0}^\infty 
\frac{(\alpha+\beta+1)_n}{(\beta+1)_n}\rho^n
P_n^{(\alpha,\beta)}(x),
\label{genGeg1}
\end{eqnarray}
where we have used the definitions
\[
\rmR=\rmR(\rho,x):=\sqrt{1+\rho^2-2\rho x},\quad
\zeta_{\pm}=\zeta_{\pm}(\rho,x):=\frac{1\pm\rho}{\sqrt{1+\rho^2-2\rho x}}.
\]
Note that the restriction given by $\rho\in\D\setminus(-1,0]$ is so that the values
of $\rho$ are ensured to remain in the domain of $P_\nu^\mu$, but may otherwise 
be relaxed to $\D$ by analytic continuation if one uses the 
Gauss hypergeometric representation.  The Ferrers function of the first kind
representation given below provides the analytic continuation to the 
segment $(-1,0]$.
Here $P_\nu^\mu:\C\setminus(-\infty,1]\to\C$ is the associated Legendre 
function of the first kind 
(see Chapter 14 in 
\cite{NIST}),
which can be defined in terms of the Gauss hypergeometric
function as follows 
\cite[(14.3.6), (15.2.2), \S 14.21(i)]{NIST}
\begin{equation}
P_\nu^\mu(z):=\frac{1}{\Gamma(1-\mu)}
\left(\frac{z+1}{z-1}\right)^{\mu/2}
\,{}_2F_1\left(
\begin{array}{c}
-\nu,\nu+1
\\[0.1cm]
1-\mu
\end{array};
\frac{1-z}{2}
\right).
\label{LegendrePdefn}
\end{equation}
The associated Legendre function of the first kind can also be 
expressed in terms of the Gauss hypergeometric function as
(see 
\cite[(14.3.18), \S 14.21(iii)]{NIST}), 
namely
\begin{equation}
P_\nu^\mu\left(z\right)=
\frac{2^\mu z^{\nu+\mu}}{\Gamma\left(1-\mu\right)\left(z^2-1\right)^{\mu/2}}
~{_2}F_1\left(
\begin{array}{c}
\frac{-\nu-\mu}{2},\frac{-\nu-\mu+1}{2}\\[0.2cm]
1-\mu
\end{array};1-\frac{1}{z^2}
\right),
\label{2F1toP}
\end{equation}
where $|\arg(z-1)|<\pi$.  We have used (\ref{2F1toP}) to re-express
the generating function (\ref{genGeg1}).
We will refer to a companion identity as one which is produced by 
applying the map $x\mapsto-x$ to an expansion over 
Jacobi, Gegenbauer, Chebyshev, or Legendre polynomials
with argument $x$,
in conjunction with the parity relations for those orthogonal polynomials.

\medskip
\noindent In our first possibility for generalizing the generating function for
Gegenbauer polynomials to a generating function for Jacobi polynomials, 
namely (\ref{genGeg1}),
we use the parity relation for Jacobi polynomials (see for instance 
\cite[Table 18.6.1]{NIST})
\begin{equation}
P_n^{(\alpha,\beta)}(-x)=(-1)^nP_n^{(\beta,\alpha)}(x),
\label{ParityJacobi}
\end{equation}
and the replacement $\alpha,\beta\mapsto\beta,\alpha$.
This produces a companion identity which is the second possibility 
for generalizing the generating function for Gegenbauer polynomials to a generating
function for Jacobi polynomials for $\rho\in(0,1)$ by
\begin{eqnarray}
&&\hspace{-1.2cm}\frac{1}{(1-\rho)^{\alpha+\beta+1}}
\,{}_2F_1
\left(
\begin{array}{c}
\frac{\alpha+\beta+1}{2},
\frac{\alpha+\beta+2}{2}\\[0.1cm]
\alpha+1
\end{array};
\frac{-2\rho(1-x)}{(1-\rho)^2}
\right)\nonumber\\[0.2cm]
&&\hspace{-0.3cm}=\left(\frac{2}{\rho(1-x)}\right)^{\alpha/2}
\frac{\Gamma(\alpha+1)}
{\rmR^{\beta+1}}
\rmP_\beta^{-\alpha}
(\zeta_-)
=\sum_{n=0}^\infty 
\frac{(\alpha+\beta+1)_n}{(\alpha+1)_n}\rho^n
P_n^{(\alpha,\beta)}(x).
\label{genGeg1companion}
\end{eqnarray}
Note that the restriction given by $\rho\in(0,1)$ is so that the values
of $\rho$ are ensured to remain in the domain of $\rmP_\nu^\mu$, but may otherwise 
be relaxed to $\D$ by analytic continuation if one uses the 
Gauss hypergeometric representation.
Here $\mathrm{P}_\nu^\mu:(-1,1)\to\C$ is the Ferrers function of the
first kind (associated Legendre function of the first kind on the cut)
through 
\cite[(14.3.1)]{NIST}, defined as
\begin{equation}
\mathrm{P}_\nu^\mu(x):=\frac{1}{\Gamma(1-\mu)}
\left(\frac{1+x}{1-x}\right)^{\mu/2}
{}_2F_1\left(
\begin{array}{c}
-\nu,\nu+1\\[0.1cm]
1-\mu
\end{array};\frac{1-x}{2}\right).
\label{FerrersPdefnGauss2F1}
\end{equation}
The Ferrers function of the first kind can also be expressed in terms
of the Gauss hypergeometric function as 
(see 
\cite[p.~167]{MOS}), namely
\begin{equation}
{\mathrm P}_\nu^\mu(x)=\frac{2^\mu x^{\nu+\mu}}{\Gamma(1-\mu)(1-x^2)^{\mu/2}}
\,{}_2F_1\left(
\begin{array}{c}
\frac{-\nu-\mu}{2},
\frac{-\nu-\mu+1}{2}\\[0.1cm]
1-\mu
\end{array};
1-\frac{1}{x^2}
\right),
\label{Ferrers1mhalfx2}
\end{equation}
for $x\in(0,1)$.  We have used
(\ref{Ferrers1mhalfx2}) to express 
(\ref{genGeg1companion}) in terms of the Ferrers function of the first kind.
One can easily see that (\ref{genGeg1}) and (\ref{genGeg1companion}) are 
generalizations of the generating function for Gegenbauer polynomials
by taking $\alpha=\beta=\nu-1/2$.
The right-hand sides easily follow using the identification 
(\ref{JactoGeg}) and the left-hand sides follow using 
\cite[(8.6.16-17)]{Abra}.  

\medskip

There exist natural extensions of
(\ref{genGeg1}),
(\ref{genGeg1companion}) in the literature
(see 
\cite[(4.3.2)]{Ismail}). The extension corresponding
to (\ref{genGeg1companion})
is given for $\rho\in(0,1)$ by 
\begin{eqnarray}
&&\hspace{-1.4cm}\frac{(\alpha+\beta+1)(1+\rho)}{(1-\rho)^{\alpha+\beta+2}}
{}_2F_1
\left(
\begin{array}{c}
\frac{\alpha+\beta+2}{2},
\frac{\alpha+\beta+3}{2}\nonumber\\[0.1cm]
\alpha+1
\end{array};
\frac{-2\rho(1-x)}{(1-\rho)^2}
\right)\nonumber\\[0.2cm]
&&\hspace{0.5cm}=\left(\frac{2}{\rho(1-x)}\right)^{\alpha/2}
\frac{(\alpha+\beta+1)(1+\rho)\Gamma(\alpha+1)}
{\rmR^{\beta+2}}
\rmP_{\beta+1}^{-\alpha}
(\zeta_-)
\nonumber\\[0.0cm]
&&\hspace{2cm}
=\sum_{n=0}^\infty 
(2n+\alpha+\beta+1)\frac{(\alpha+\beta+1)_n}{(\alpha+1)_n}\rho^n
P_n^{(\alpha,\beta)}(x),
\label{extnJacgegnfac1}
\end{eqnarray}
and its companion identity corresponding to 
(\ref{genGeg1}), for $\rho\in\D\setminus(-1,0]$ is 
\begin{eqnarray}
&&\hspace{-0.9cm}\frac{(\alpha+\beta+1)(1-\rho)}{(1+\rho)^{\alpha+\beta+2}}
{}_2F_1
\left(
\begin{array}{c}
\frac{\alpha+\beta+2}{2},
\frac{\alpha+\beta+3}{2}\\[0.1cm]
\beta+1
\end{array};
\frac{2\rho(1+x)}{(1+\rho)^2}
\right)\nonumber\\[0.1cm]
&&\hspace{0.5cm}=\left(\frac{2}{\rho(1+x)}\right)^{\beta/2}
\frac{(\alpha+\beta+1)(1-\rho)\Gamma(\beta+1)}
{\rmR^{\alpha+2}}
P_{\alpha+1}^{-\beta}
(\zeta_+)
\nonumber\\[0.0cm]
\nonumber\\[0.0cm]
&&\hspace{2.0cm}
=\sum_{n=0}^\infty
(2n+\alpha+\beta+1)\frac{(\alpha+\beta+1)_n}{(\beta+1)_n}\rho^n
P_n^{(\alpha,\beta)}(x).
\label{extnJacgegnfac2}
\end{eqnarray}
We have used 
(\ref{2F1toP}),
(\ref{Ferrers1mhalfx2})
to re-express
these Gauss hypergeometric function generating functions as
associated Legendre functions.
We have not seen the companion identity (\ref{extnJacgegnfac2}) in
the literature, but it is an obvious consequence of \cite[(4.3.2)]{Ismail} 
using parity.  On the other hand, we have not seen the associated Legendre function 
representations 
of 
(\ref{extnJacgegnfac1}), (\ref{extnJacgegnfac2}) in the literature.

Upon examination of these two sets of generating functions, we suspected
that these were just two examples of an infinite sequence of such expansions.
This led us to the proof of the following theorem, which is a Jacobi polynomial expansion
which generalizes the generating function for Gegenbauer polynomials
(\ref{generatingfunctionforGeg}).
According to Ismail (2005)
\cite[(4.3.2)]{Ismail}, 
the generating functions 
(\ref{extnJacgegnfac1}),
(\ref{extnJacgegnfac2}), their generalizations
Theorem \ref{firstbigthm}, Corollary \ref{firstbigcor} and their 
corresponding definite integrals
(\ref{defintJac1}),
(\ref{secondintegralcor}),
are closely related to the Poisson
kernel for Jacobi polynomials, so our new generalizations
will have corresponding applications.

\begin{thm}
Let $m\in\N_0$, $\alpha,\beta>-1$ such that if $\alpha,\beta\in(-1,0)$ then
$\alpha+\beta+1\ne 0$, $x\in[-1,1],$ 
$\rho\in\D\setminus(-1,0].$ Then
\begin{equation}
\frac{(1+x)^{-\beta/2}}
{\rmR^{\alpha+m+1}}
P_{\alpha+m}^{-\beta}
(\zeta_+)
=\sum_{n=0}^\infty
a_{n,m}^{(\alpha,\beta)}(\rho)P_n^{(\alpha,\beta)}(x),
\label{generalizedexpansion}
\end{equation}
where $a_{n,m}^{(\alpha,\beta)}:
\D\setminus(-1,0]\to\C$ is defined by
\begin{eqnarray*}
&&a_{n,m}^{(\alpha,\beta)}(\rho):=
\frac{(2n+\alpha+\beta+1) \Gamma(\alpha+\beta+n+1) (\alpha+\beta+m+1)_{2n}}
{ 
2^{\beta/2} 
\Gamma(\beta+n+1)
}\\[0.2cm]
&&\hspace{5cm}\times\frac{1}
{\rho^{(\alpha+1)/2} 
(1-\rho)^m}
P_{-m}^{-\alpha-\beta-2n-1}
\left(\frac{1+\rho}{1-\rho}\right).
\end{eqnarray*}
\label{firstbigthm}
\end{thm}
\noindent {\bf Proof.} 
Let 
$\rho\in(0,\epsilon)$ with $\epsilon$ sufficiently small.
Then using the definition of the following Gauss 
hypergeometric function 
\begin{eqnarray}
&&\hspace{-1.1cm}{}_2F_1\left(
\begin{array}{c}
\frac{\alpha+\beta+m+1}{2},\frac{\alpha+\beta+m+2}{2}
\\[0.1cm]
\beta+1
\end{array}
;\frac{2\rho(1+x)}{(1+\rho)^2}
\right)\nonumber\\[0.2cm]
&&\hspace{2cm}=\sum_{n=0}^\infty\frac
{
\left(\frac{\alpha+\beta+m+1}{2}\right)_n
\left(\frac{\alpha+\beta+m+2}{2}\right)_n(2\rho)^n(1+x)^n
}{n!(\alpha+1)_n(1+\rho)^{2n}},
\label{gaussstart}
\end{eqnarray}
the expansion of $(1+x)^n$ in terms of Jacobi polynomials
is given by
\begin{equation}
(1+x)^n=2^n(\beta+1)_n\sum_{k=0}^n
\frac{(-1)^k(-n)_k(\alpha+\beta+2k+1)(\alpha+\beta+1)_k}
{(\alpha+\beta+1)_{n+k+1}(\beta+1)_k}
P_k^{(\alpha,\beta)}(x),
\label{1mxnJacobi}
\end{equation}
whose coefficients can be determined using 
orthogonality of Jacobi polynomials 
(see Appendix)
combined with the Mellin transform given in 
\cite[(18.17.36)]{NIST}.
By inserting (\ref{1mxnJacobi}) in the right-hand side of (\ref{gaussstart}), we obtain
an expansion of the Gauss hypergeometric function on the left-hand side of (\ref{gaussstart})
in terms of Jacobi polynomials.  By interchanging the two sums (with justification by
absolute convergence), shifting the $n$-index by $k$, and taking advantage of standard 
properties such as 
\begin{eqnarray*}
(-n-k)_k=\frac{(-1)^k(n+k)!}{n!},\\[0.2cm]
(a)_{n+k}=(a)_k (a+k)_n, \\[0.2cm]
\left(\frac{a}{2}\right)_n\left(\frac{a+1}{2}\right)_n=\frac{1}{2^{2n}}\left(a\right)_{2n},
\end{eqnarray*}
$n,k\in\N_0$, $a\in\C$, produces a Gauss hypergeometric function as the coefficient of the 
Jacobi polynomial expansion.  The resulting expansion is
\begin{equation}
{}_2F_1
\left(
\begin{array}{c}
\frac{\alpha+\beta+m+1}{2},
\frac{\alpha+\beta+m+2}{2}\\[0.1cm]
\beta+1
\end{array};
\frac{2\rho(1+x)}{(1+\rho)^2}
\right)=\sum_{n=0}^\infty f_{n,m}^{(\alpha,\beta)}(\rho)P_n^{(\alpha,\beta)}(x),
\label{generalizedexpansiongausshyper}
\end{equation}
where 
$f_{n,m}^{(\alpha,\beta)}:(0,\epsilon)\to\R$ is defined by
\begin{eqnarray*}
&&
f_{n,m}^{(\alpha,\beta)}(\rho):=
\frac
{
(2n+\alpha+\beta+1)(\alpha+\beta+1)_n(\alpha+\beta+m+1)_{2n}\,\rho^n
}
{(\beta+1)_n
(\alpha+\beta+1)_{2k+1}(1+\rho)^{2n}
}\\[0.2cm]
&&\hspace{3.0cm}\times
\,{}_2F_1\left(
\begin{array}{c}
\frac{\alpha+\beta+m+2n+1}{2},\frac{\alpha+\beta+m+2n+2}{2}\\[0.1cm]\alpha+\beta+2n+2
\end{array};
\frac{4\rho}{(1+\rho)^2}
\right).
\end{eqnarray*}
The above expansion is actually analytic on $\D$.  However, if we express it in
terms of associated Legendre functions, then we must necessarily subdivide it
into two regions.
The Gauss hypergeometric function coefficient of this expansion, as well as the 
Gauss hypergeometric function on the left-hand of 
(\ref{generalizedexpansiongausshyper}) are realized to be associated Legendre functions 
of the first kind through 
(\ref{2F1toP}).
Both sides of the resulting Jacobi polynomial expansion 
are analytic
functions on $\rho\in\D\setminus(-1,0].$ 
Since we know that (\ref{generalizedexpansion}) is valid 
for $\rho\in(0,\epsilon),$ then by the identity theorem for 
analytic functions, the equation holds on this domain.  This completes 
the proof.  $\hfill\blacksquare$

\medskip

\noindent Note that the left-hand side of Theorem \ref{firstbigthm} can be rewritten as
\[
\frac{(\rho/2)^{\beta/2}}{\Gamma(\beta+1)(1+\rho)^{\alpha+\beta+m+1}}\,
{}_2F_1\left(
\begin{array}{c}
\frac{\alpha+\beta+m+1}{2},\frac{\alpha+\beta+m+2}{2}\\[0.2cm]\beta+1
\end{array};
\frac{2\rho(1+x)}{(1+\rho)^2}\right).
\]
\medskip

\noindent We have also derived the companion identity to (\ref{generalizedexpansion}), which
we give in the following corollary.

\begin{cor}
Let $m\in\N_0$, $\alpha,\beta>-1$ such that if $\alpha,\beta\in(-1,0)$ then
$\alpha+\beta+1\ne 0$, $x\in[-1,1]$, 
$\rho\in(0,1)$. Then
\begin{equation}
\frac{(1-x)^{-\alpha/2}}
{\rmR^{\beta+m+1}}
{\mathrm P}_{\beta+m}^{-\alpha}
(\zeta_-)
=
\sum_{n=0}^\infty b_{n,m}^{(\alpha,\beta)}(\rho)P_n^{(\alpha,\beta)}(x),
\label{companiongeneJacobi}
\end{equation}
where $b_{n,m}^{(\alpha,\beta)}:(0,1)\to\R$ is defined by
\begin{eqnarray*}
b_{n,m}^{(\alpha,\beta)}(\rho):=
\frac{(2n+\alpha+\beta+1) \Gamma(\alpha+\beta+n+1) (\alpha+\beta+m+1)_{2n}}
{ 
2^{\alpha/2} 
\Gamma(\alpha+n+1)}
\nonumber\\[0.2cm]
\times\frac{1}
{
\rho^{(\beta+1)/2}
(1+\rho)^m }
\mathrm{P}_{-m}^{-\alpha-\beta-2n-1}
\left(\frac{1-\rho}{1+\rho}\right).
\end{eqnarray*}
\label{firstbigcor}
\end{cor}
\noindent {\bf Proof.} 
We start with (\ref{generalizedexpansiongausshyper}) and apply the
parity relation for Jacobi polynomials
(\ref{ParityJacobi}).  Let $\rho\in(0,1)$. 
The Gauss hypergeometric function 
coefficient of the 
Jacobi expansion is seen to be a Ferrers function of the first kind 
(\ref{FerrersPdefnGauss2F1}).
After the application of the parity relation, the left-hand side also reduces 
to a Ferrers function of the first kind 
through (\ref{Ferrers1mhalfx2}).
This completes the proof.
$\hfill\blacksquare$

\bigskip
\noindent Theorem 
\ref{firstbigthm}
generalizes (\ref{genGeg1}), 
(\ref{extnJacgegnfac2}),
while Corollary \ref{firstbigcor}
generalizes
(\ref{genGeg1companion}), 
(\ref{extnJacgegnfac1}). 
Both Theorem \ref{firstbigthm} and Corollary \ref{firstbigcor}
generalize the generating function for Gegenbauer polynomials
(\ref{generatingfunctionforGeg}), which is its own companion identity.  

\medskip

\subsection{Expansions and definite integrals from the Szeg\H{o} transformation}
\label{ExpansionsanddefiniteintegralsfromtheSzegotransformation}

If one applies on the complex plane, the Szeg\H{o} transformation (conformal map)
\begin{equation}
z=\frac{1+\rho^2}{2\rho},
\label{confmap}
\end{equation}
(which maps a circle with radius less than unity to an ellipse with the foci at $\pm 1$)
to the expansion in Theorem \ref{firstbigthm}, then one obtains a new expansion.
By 
\cite[Theorem 12.7.3]{Szego}, this new Jacobi polynomial 
expansion is convergent for all $x\in\C$ within the interior of this ellipse.
Applying (\ref{confmap}) to (\ref{generalizedexpansion}) yields the following 
corollary.

\begin{cor}
Let
$m\in\N_0$,
$\alpha,\beta>-1$ such that if $\alpha,\beta\in(-1,0)$ then $\alpha+\beta+1\ne 0$,
$x,z\in\C$, with $z\in\C\setminus(-\infty,1]$ on any ellipse with the foci 
at $\pm 1$ and $x$ in the interior of that ellipse. Then
\begin{equation}
\frac{(1+x)^{-\beta/2}}{(z-x)^{(\alpha+m+1)/2}}
P_{\alpha+m}^{-\beta}\left(
\frac{1+z-\sqrt{z^2-1}}
{\sqrt{2(z-\sqrt{z^2-1})(z-x)}}
\right)=\sum_{n=0}^\infty c_{n,m}^{(\alpha,\beta)}(z)
P_n^{(\alpha,\beta)}(x),
\label{Szegocorollaryeqn}
\end{equation}
where $c_{n,m}^{(\alpha,\beta)}:\C\setminus(-\infty,1]\to\C$  
is defined by
\begin{eqnarray*}
&&c_{n,m}^{(\alpha,\beta)}(z):=
\frac{(2n+\alpha+\beta+1)\Gamma(\alpha+\beta+n+1)(\alpha+\beta+m+1)_{2n}}
{
2^{(\beta-\alpha-m-1)/2}
\Gamma(\beta+n+1)}\\[0.2cm]
&&
\hspace{4cm}\times 
\frac{(z-\sqrt{z^2-1})^{m/2}}{\left(1-z+\sqrt{z^2-1}\right)^m}
P_{-m}^{-\alpha-\beta-2n-1}
\left(\sqrt{\frac{z+1}{z-1}}\right).
\end{eqnarray*}
\label{Szegocorollary}
\end{cor}


We would just like to briefly note that one may use the Szeg\H{o} transformation 
(\ref{confmap})
to obtain new expansion formulae and corresponding definite integrals 
from all the Jacobi, Gegenbauer, Legendre and 
Chebyshev polynomial expansions used in this paper.  For the sake of brevity,
we leave this to the reader.

\section{Expansions over Gegenbauer polynomials}
\label{Gegenbauerpolynomials}

The Gegenbauer polynomials $C_n^{\mu}:\C\to\C$ can be defined 
in terms of the terminating Gauss hypergeometric series as follows
(\cite[(18.5.9)]{NIST})
\begin{equation}
C_n^{\mu}(z):=\frac{(2\mu)_n}{n!}\,{}_2F_1\left(
\begin{array}{c}
-n,n+2\mu\\[0.1cm]
\mu+\frac12
\end{array};
\frac{1-z}{2}
\right),
\label{gegenbauerdefngauss2F1}
\end{equation}
for $n\in\N_0$ and 
$\mu\in(-1/2,\infty)\setminus\{0\}.$
The Gegenbauer polynomials (\ref{gegenbauerdefngauss2F1}) are defined 
for $\mu\in(-1/2,\infty)\setminus\{0\}.$
However many of the formulae listed below actually make sense in the limit
as $\mu\to0$.
In this case, one should take the limit of the expression as $\mu\to0$ 
with the interpretation of obtaining Chebyshev polynomials of the first kind
(see \S\ref{Chebyshevpolynomialsofthefirstkind} for the details of this 
limiting procedure).

\begin{cor}
Let $m\in\N_0,$ 
$\mu\in(-1/2,\infty)\setminus\{0\},$
$x\in[-1,1]$. 
If $\rho\in\D\setminus(-1,0],$ then 
\begin{equation}
\frac{1}
{\rmR^{2\mu+m}}
C_m^\mu
(\zeta_+)
=\frac{2\Gamma(2\mu+m)}{m!\rho^\mu(1-\rho)^m}
\sum_{n=0}^\infty(n+\mu)(2\mu+m)_{2n}P_{-m}^{-2\mu-2n}\left(\frac{1+\rho}{1-\rho}\right)
C_n^\mu(x),
\label{gegenposintegral}
\end{equation}
and if $\rho\in(0,1)$ then
\begin{equation}
\frac{1}
{\rmR^{2\mu+m}}
C_m^\mu
(\zeta_-)
=\frac{2\Gamma(2\mu+m)}{m!\rho^\mu(1+\rho)^m}
\sum_{n=0}^\infty(n+\mu)(2\mu+m)_{2n}{\mathrm P}_{-m}^{-2\mu-2n}\left(\frac{1-\rho}{1+\rho}\right)
C_n^\mu(x).
\label{negdefintGegenabuer}
\end{equation}
\end{cor}
\noindent {\bf Proof.} 
Using (\ref{generalizedexpansion}), substitute $\alpha=\beta=\mu-1/2$ along with 
(\ref{JactoGeg}) and
\cite[(14.3.22)]{NIST}, namely
\[
P_{n+\mu-1/2}^{1/2-\mu}(z)=
\frac{2^{\mu-1/2}\Gamma(\mu)n!}
{\sqrt{\pi}\,\Gamma(2\mu+n)}
(z^2-1)^{\mu/2-1/4}\,C_n^\mu(z).
\]
Through (\ref{LegendrePdefn}), we see that the Gauss hypergeometric 
function in the definition of the associated Legendre function of the
first kind on the right-hand side is terminating and therefore defines 
an analytic function
for $\rho\in\D$.  The analytic continuation to the segment $\rho\in(0,1]$
is provided by replacing the associated Legendre function of the first
kind with the Ferrers function of the first kind with argument
$(1-\rho)/(1+\rho)$.
$\hfill\blacksquare$

%

\medskip
As an example for re-expression using the elementarity of associated Legendre functions
which was mentioned in the introduction, we now apply to two generating function
results of Koekoek {\it et al.} (2010) \cite[(9.8.32)]{Koekoeketal} 
and Rainville (1960) \cite[(144.8)]{Rainville}.

\begin{thm}
Let $\lambda\in\C,$ 
$\mu\in(-1/2,\infty)\setminus\{0\},$
$\rho\in(0,1),$ $x\in[-1,1]$. Then
\begin{eqnarray}
&&\hspace{-0.9cm}(1-x^2)^{1/4-\mu/2}
P_{\mu-\lambda-1/2}^{1/2-\mu}
\left(\rmR+\rho\right)\,
{\mathrm P}_{\mu-\lambda-1/2}^{1/2-\mu}
\left(
\rmR
-\rho\right)\nonumber\\[0.2cm]
&&\hspace{+2.5cm}=\frac{2^{1/2-\mu}}{\Gamma(\mu+\frac12)}
\sum_{n=0}^\infty \frac{(\lambda)_n\,(2\mu-\lambda)_n}
{(2\mu)_n\,\Gamma(\mu+\frac12+n)}\rho^{\mu-1/2+n}C_n^\mu(x).
\label{GegenforInt2}
\end{eqnarray}
\label{Koekoekgenfun}
\end{thm}
\noindent {\bf Proof.}
\cite[(9.8.32)]{Koekoeketal} give a generating function for Gegenbauer polynomials, namely
\begin{eqnarray}
&&\hspace{-0.8cm}{_2}F_1\left(
\begin{array}{c}
\lambda,2\mu-\lambda\\[0.2cm]
\mu+\frac12
\end{array};\frac{1-\rmR
-\rho}{2}\right)
{_2}F_1\left(
\begin{array}{c}
\lambda,2\mu-\lambda\\[0.2cm]
\mu+\frac12
\end{array};\frac{1-\rmR+\rho}{2}\right)\nonumber\\[0.2cm]
&&\hspace{+5.7cm}=\sum_{n=0}^\infty \frac{(\lambda)_n\,(2\mu-\lambda)_n}
{(2\mu)_n\,(\mu+\frac12)_n}\rho^nC_n^\mu(x).\nonumber
\end{eqnarray}
Using 
\cite[(15.8.17)]{NIST} to do a quadratic transformation
on the Gauss hypergeometric functions and then using
(\ref{2F1toP}), with degree and order given by $\mu-\lambda-1/2,$ $1/2-\mu,$ 
respectively and 
either $z=\rmR+\rho$ or $z=\rmR-\rho$.
Simplification completes the proof.
$\hfill\blacksquare$

\begin{thm}
Let $\alpha\in\C,$ 
$\mu\in(-1/2,\infty)\setminus\{0\},$
$\rho\in(0,1),$ $x\in[-1,1]$. Then
\begin{eqnarray}
&&\hspace{-0.9cm}\frac{\left(1-x^2\right)^{1/4-\mu/2}}
{\rmR^{1/2+\alpha-\mu}}
{\mathrm P}_{\mu-\alpha-1/2}^{1/2-\mu}\left(\frac{1-\rho x}
{\rmR}
\right)
=\frac{(\rho/2)^{\mu-1/2}}{\Gamma(\mu+\frac12)}\sum_{n=0}^\infty  \frac{(\alpha)_n}{(2\mu)_n}\rho^nC_n^\mu(x).
\label{GegenforInt3}
\end{eqnarray}
\end{thm}
\noindent {\bf Proof.}
On p.~279 of 
\cite[(144.8)]{Rainville} there is a generating function
for Gegenbauer polynomials, namely
\[
(1-\rho x)^{-\alpha}\, {_2}F_1\left(
\begin{array}{c}
\frac\alpha2,\frac{\alpha+1}{2}\\[0.2cm]
\mu+\frac12
\end{array};\frac{-\rho^2(1-x^2)}{(1-\rho x)^2}\right)
=\sum_{n=0}^\infty\frac{(\alpha)_n}{(2\mu)_n}\rho^nC_n^\mu(x).\nonumber
\]
Using (\ref{2F1toP}) to rewrite the Gauss hypergeometric function on the left-hand side of the above equation completes the proof.
$\hfill\blacksquare$

\section{Expansions over Chebyshev polynomials of the second kind}
\label{Chebyshevpolynomialsofthesecondkind}

The Chebyshev polynomials of the second kind can be obtained from the Gegenbauer
polynomials using 
\cite[(18.7.4)]{NIST}, namely
\begin{equation}
U_n(z)=C_n^1(z),
\label{UtoC}
\end{equation}
for $n\in\N_0$.
Hence and through (\ref{gegenbauerdefngauss2F1}),
the Chebyshev polynomials of the second kind $U_n:\C\to\C$ can be defined 
in terms of the terminating Gauss hypergeometric series as follows
\begin{equation}
U_n(z):=(n+1)\,{}_2F_1\left(\begin{array}{c}
-n,n+2\\[0.1cm]
\frac32
\end{array};
\frac{1-z}{2}
\right).
\label{Chebyshev2defnGauss2F1}
\end{equation}

\begin{cor} 
Let $m\in\N_0,$ 
$x\in[-1,1].$
If $\rho\in\D\setminus(-1,0]$ then
\begin{eqnarray}
&&\hspace{-1.3cm}\frac{1}
{\rmR^{m+2}}
U_m
(\zeta_+)
=\frac{2(m+1)}{\rho(1-\rho)^m}
\sum_{n=0}^\infty(n+1)(m+2)_{2n}P_{-m}^{-2n-2}\left(\frac{1+\rho}{1-\rho}\right)
U_n(x),
\label{posexpansioncheby2}
\end{eqnarray}
and if $\rho\in(0,1)$ then
\begin{eqnarray}
&&\hspace{-1.3cm}\frac{1}
{\rmR^{m+2}}
U_m
(\zeta_-)
=\frac{2(m+1)}{\rho(1+\rho)^m}
\sum_{n=0}^\infty(n+1)(m+2)_{2n}{\mathrm P}_{-m}^{-2n-2}\left(\frac{1-\rho}{1+\rho}\right)
U_n(x).
\end{eqnarray}
\end{cor}
\noindent {\bf Proof.} 
Using (\ref{gegenposintegral}), with (\ref{UtoC}) and
\begin{equation}
U_m(z)=\sqrt{\frac{\pi}{2}}\frac{m+1}{(z^2-1)^{1/4}}P_{m+1/2}^{-1/2}(z),
\label{Cheby2defnLegP}
\end{equation}
which follows from 
(\ref{LegendrePdefn}), 
(\ref{Chebyshev2defnGauss2F1}),
and \cite[(15.8.1)]{NIST}.
Analytically continuing to the segment $\rho\in(0,1)$ completes 
the proof. $\hfill\blacksquare$

\medskip

\noindent Note that using 
\cite[(8.6.9)]{Abra}, namely
\[
P_\nu^{-1/2}\left(z\right)=
\sqrt{\frac{2}{\pi}}\frac{\left(z^2-1\right)^{-1/4}}{\left(2\nu+1\right)}
\left[\left(z+\sqrt{z^2-1}\right)^{\nu+1/2}-\left(z+\sqrt{z^2-1}\right)^{-\nu-1/2}\right],
\]
and
(\ref{Cheby2defnLegP}) one can derive the 
elementary function representation for the 
Chebyshev polynomials of the second kind
\cite[(1.52)]{MasonHandscomb}.


\begin{cor} Let $\rho\in\D,$ $x\in[-1,1]$. Then
\begin{eqnarray}
&&\hspace{-1.0cm}(1-x^2)^{-1/4} 
P_{1/2-\lambda}^{-1/2}
\left(
\rmR
+\rho\right)
{\mathrm P}_{1/2-\lambda}^{-1/2}\left(
\rmR
-\rho\right)\nonumber\\[0.2cm]
&&\hspace{+3.1cm}=\frac{2^{5/2}\sqrt{\rho}}{\pi}
\sum_{n=0}^\infty \frac{(\lambda)_n\,(2-\lambda)_n\,2^{2n}\rho^n}{(2n+2)!}U_n(x).
\label{ConnorUforint1}
\end{eqnarray}
\end{cor}
\noindent {\bf Proof.}
Substituting $\mu=1$ into (\ref{GegenforInt2}),
and using (\ref{UtoC}) with simplification, completes the proof.
$\hfill\blacksquare$

\begin{cor}
Let $\alpha\in\C,$ $\rho\in(0,1),$ $x\in[-1,1]$. Then
\begin{equation}
\frac{
\rmR^{1/2-\alpha}
}
{(1-x^2)^{1/4}}
{\mathrm P}_{1/2-\alpha}^{-1/2}\left(\frac{1-\rho x}
{\rmR}
\right)
=\sqrt{\frac{2\rho}{\pi}}\sum_{n=0}^\infty\frac{(\alpha)_n}{(n+1)!}\rho^nU_n(x).
\label{ConnorUforint4}
\end{equation}
\end{cor}
\noindent {\bf Proof.}
Using (\ref{UtoC}), and substituting $\mu=1$ into (\ref{GegenforInt3}) with 
simplification,
produces this generating function for Chebyshev polynomials of the second kind.
$\hfill\blacksquare$

\section{Expansions over Legendre polynomials}
\label{Legendrepolynomialsofthefirstkind}

Legendre polynomials can be obtained from the Gegenbauer
polynomials using 
\cite[(18.7.9)]{NIST}, namely
\begin{equation}
P_n(z)=C_n^{1/2}(z),
\label{PtoC}
\end{equation}
for $n\in\N_0$.
Hence and through (\ref{gegenbauerdefngauss2F1}),
the Legendre polynomials of the second kind $P_n:\C\to\C$ can be defined 
in terms of the terminating Gauss hypergeometric series as follows
\begin{equation}
P_n(z):={}_2F_1\left(\begin{array}{c}
-n,n+1\\[0.1cm]
1
\end{array};
\frac{1-z}{2}
\right).
\label{defnLegendrepolygauss2F1}
\end{equation}
Using (\ref{PtoC}) we can write the previous expansions over Gegenbauer polynomials in terms
of expansions over Legendre polynomials.

\begin{cor}
Let $m\in\N_0,$ 
$x\in[-1,1]$. 
If $\rho\in\D\setminus(-1,0]$ then
\begin{eqnarray}
&&\hspace{-1.2cm}\frac{1}
{\rmR^{m+1}}
P_m
(\zeta_+)
=\frac{(1-\rho)^{-m}}{\sqrt{\rho}}
\sum_{n=0}^\infty(2n+1)(m+1)_{2n}P_{-m}^{-2n-1}\left(\frac{1+\rho}{1-\rho}\right)
P_n(x),
\label{Legendrepolyposexpansions}
\end{eqnarray}
and if $\rho\in(0,1)$ then
\begin{eqnarray}
&&\hspace{-1.2cm}\frac{1}
{\rmR^{m+1}}
P_m
(\zeta_-)
=\frac{(1+\rho)^{-m}}{\sqrt{\rho}}
\sum_{n=0}^\infty(2n+1)(m+1)_{2n}{\mathrm P}_{-m}^{-2n-1}\left(\frac{1-\rho}{1+\rho}\right)
P_n(x).
\label{Legendrepolynegexpansions}
\end{eqnarray}
\end{cor}
\noindent {\bf Proof.} 
Using (\ref{gegenposintegral}), substitute $\mu=1/2$ with 
(\ref{PtoC}) 
and
$P_m(z)=P_m^0(z),$
which follows from 
(\ref{LegendrePdefn}),
(\ref{defnLegendrepolygauss2F1}).  Analytic continuation to 
$\rho\in(0,1)$ completes the proof. 
$\hfill\blacksquare$


\begin{cor}
Let $\lambda\in\C,$ $\rho\in\{x\in\C:|z|<1\},$ $x\in[-1,1]$. Then
\begin{equation}
P_{-\lambda}
\left(\rmR+\rho\right)
{\mathrm P}_{-\lambda}
\left(\rmR-\rho\right)
=\sum_{n=0}^\infty\frac{(\lambda)_n\,(1-\lambda)_n}{(n!)^2}\rho^nP_n(x).
\label{LegPforInt1}
\end{equation}
\label{Braftheoremhere}
\end{cor}
\noindent {\bf Proof.}
Substituting $\mu=1/2$ into (\ref{GegenforInt2}) and using (\ref{PtoC})
with simplification completes the proof.  $\hfill\blacksquare$
 
\medskip

\noindent Note that Corollary \ref{Braftheoremhere} is just a restatement of 
\cite[Theorem A]{WanWad12}, and therefore 
Theorem \ref{Koekoekgenfun} is a generalization of Brafman's theorem.

\begin{cor}
Let $\alpha\in\C,$ $\rho\in\D,$ $x\in[-1,1].$ Then
\begin{equation}
\rmR^{-\alpha}
\,{\mathrm P}_{\alpha-1}\left(\frac{1-\rho x}
{\rmR}
\right)
=\sum_{n=0}^\infty \frac{(\alpha)_n}{n!}\rho^nP_n(x). 
\label{LegPforInt2}
\end{equation}
\end{cor}
\noindent {\bf Proof.}
Substituting $\mu=1/2$ in (\ref{GegenforInt3}) with simplification completes the proof.
$\hfill\blacksquare$

\medskip
As a further example of the elementarity of associated Legendre functions mentioned
in the introduction, we apply to the recent generating function results 
of Wan \& Zudelin (2012) \cite{WanWad12}.

\begin{thm}
Let $x,y$ be in a neighborhood of 1.  Then
\begin{eqnarray*}
&&\hspace{0.1cm}\frac{\pi^2}{2}
\sum_{n=0}^\infty\frac{(\frac12)_n^2}{(n!)^2}
P_{2n}\left(\frac{(x+y)(1-xy)}{(x-y)(1+xy)}\right)
\left(\frac{x-y}{1+xy}\right)^{2n}
\nonumber\\[0.2cm]
&&\hspace{1cm}
=\left\{
\begin{array}{ll}
\displaystyle \frac{\pi^2}{2}
& \mathrm{if}\ x=y=1,\\[0.4cm]
\displaystyle \frac{1+xy}{xy}\,
K\left(\frac{\sqrt{x^2-1}}{x}\right)
K\left(\frac{\sqrt{y^2-1}}{y}\right)
& \mathrm{if}\ x,y\ge 1,\\[0.6cm]
\displaystyle \frac{1+xy}{x}\,
K\left(\frac{\sqrt{x^2-1}}{x}\right)
K\left(\sqrt{1-y^2}\right)
& \mathrm{if}\ x\ge 1
{\rm \ and\ } 
y\le 1,\\[0.6cm]
\displaystyle \frac{1+xy}{y}\,
K\left(\sqrt{1-x^2}\right)
K\left(\frac{\sqrt{y^2-1}}{y}\right)
& \mathrm{if}\ x\le 1
{\rm \ and\ } 
y\ge 1,\\[0.6cm]
\displaystyle (1+xy)\,
K\left(\sqrt{1-x^2}\right)
K\left(\sqrt{1-y^2}\right)
& \mathrm{if}\ x,y\le 1.
\end{array}
\right. 
\end{eqnarray*}
\label{Braf1}
\end{thm}
\noindent {\bf Proof.} 
If we start with (10) from 
\cite{WanWad12}, namely 
\begin{eqnarray*}
&&\sum_{n=0}^\infty\frac{\left(\frac12\right)_n^2}{(n!)^2}
P_{2n}\left(\frac{(x+y)(1-xy)}{(x-y)(1+xy)}\right)
\left(\frac{x-y}{1+xy}\right)^{2n}\nonumber\\[0.2cm]
&&\hspace{4.0cm}=\frac{1+xy}{2}
\,{}_2F_1\left(\begin{array}{c}\frac12,\frac12\\[0.2cm]1\end{array};1-x^2\right)
{}_2F_1\left(\begin{array}{c}\frac12,\frac12\\[0.2cm]1\end{array};1-y^2\right),
\end{eqnarray*}
and use \cite[(15.9.21)]{NIST} we can express the Gauss hypergeometric functions
as Legendre functions. For instance
\[
{}_2F_1\left(\begin{array}{c}\frac12,\frac12\\[0.2cm]1\end{array};1-x^2\right)
=P_{-1/2}(2x^2-1),
\]
with $x\in\C\setminus(-\infty,0]$.
This domain is because the Legendre function of the first kind $P_\nu$ and the Ferrers function of 
the first kind $\rmP_\nu$, both with order $\mu=0,$ are given by the same Gauss 
hypergeometric function and are continuous across argument unity
(cf.~(\ref{LegendrePdefn}), (\ref{FerrersPdefnGauss2F1})).  So there 
is no distinction between these two functions, except that the Ferrers function
has argument on the real line with modulus less than unity and the 
Legendre function is defined on $\C\setminus(-\infty,1)$ (both being well defined
with argument unity).  (Hence there really is no need to use two different
symbols to denote this function.)  The proof is completed by noting the two formulae
\[
P_{-1/2}(z)=\frac{2}{\pi}\sqrt{\frac{2}{z+1}}K\left(\sqrt{\frac{z-1}{z+1}}\right),
\]
\[
\rmP_{-1/2}(x)=\frac{2}{\pi}K\left(\sqrt{\frac{1-x}{2}}\right)
\]
\cite[(8.13.1), (8.13.8)]{Abra}, where
$K:[0,1)\to[\pi/2,\infty)$ is the complete elliptic integral of the first
kind defined by \cite[(19.2.8)]{NIST}
\[
K(k):=\frac{\pi}{2}\,{}_2F_1\left(\begin{array}{c}\frac12,\frac12\\[0.2cm]1\end{array};k^2\right).
\]
$\hfill\blacksquare$

\medskip

\begin{thm}
Let $x,y$ be in a neighborhood of 1.  Then
\begin{eqnarray*}
&&=3\sum_{n=0}^\infty
\frac{\left(\frac13\right)_n\left(\frac23\right)_n}{(n!)^2}
P_{3n}\left(\frac{x+y-2x^2y^2}{(x-y)\sqrt{1+4xy(x+y)}}\right)
\left(\frac{x-y}{\sqrt{1+4xy(x+y)}}\right)^{3n}.\nonumber\\[0.4cm]
&&=
\sqrt{1+4xy(x+y)}
\left\{\begin{array}{ll}
1
& \mathrm{if}\ x=y=1 \\[0.2cm]
P_{-1/3}\left(2x^3-1\right)P_{-1/3}\left(2y^3-1\right) 
& \mathrm{if}\ x,y> 1  \\[0.2cm]
\rmP_{-1/3}\left(2x^3-1\right)P_{-1/3}\left(2y^3-1\right) 
& \mathrm{if}\ x< 1 {\rm \ and\ } y > 1 \\[0.2cm]
P_{-1/3}\left(2x^3-1\right)\rmP_{-1/3}\left(2y^3-1\right) 
& \mathrm{if}\ x>1 {\rm \ and\ } y<1 \\[0.2cm]
\rmP_{-1/3}\left(2x^3-1\right)\rmP_{-1/3}\left(2y^3-1\right) 
& \mathrm{if}\ x,y < 1 
\end{array}
\right.
\end{eqnarray*}
\end{thm}
\noindent {\bf Proof.}  This follows by 
\cite[(11)]{WanWad12} and \cite[(15.9.21)]{NIST}. 
$\hfill\blacksquare$

\section{Expansions over Chebyshev polynomials of the first kind}
\label{Chebyshevpolynomialsofthefirstkind}

The Chebyshev polynomials of the first kind $T_n:\C\to\C$ can be defined 
in terms of the terminating Gauss hypergeometric series as follows
(
\cite[p.~257]{MOS})
\begin{equation}
T_n(z):={}_2F_1\left(\begin{array}{c}
-n,n\\[0.1cm]
\frac12
\end{array};
\frac{1-z}{2}
\right),
\label{defnChebyshev1polygauss2F1}
\end{equation}
for $n\in\N_0$.
The Chebyshev polynomials of the first kind can be obtained from the Gegenbauer
polynomials using 
\cite[(6.4.13)]{AAR}, namely
\begin{equation}
T_n(z)=\frac{1}{\epsilon_n}\lim_{\mu\to 0}\frac{n+\mu}{\mu}C_n^\mu(z),
\label{defnChebyshev1intermsofGegenbauer}
\end{equation}
where the Neumann factor $\epsilon_n\in\{1,2\}$, commonly seen in Fourier 
cosine series, is defined as 
$\epsilon_n:=2-\delta_{n,0}.$

\begin{cor}
Let $m\in\N_0,$ 
$x\in[-1,1]$. 
If $\rho\in\D\setminus(-1,0]$ then
\begin{equation}
\hspace{-0.4cm}\frac{1}
{\rmR^{m}}
T_m
(\zeta_+)
=\frac{1}{(1-\rho)^m}
\sum_{n=0}^\infty\epsilon_n(m)_{2n}P_{-m}^{-2n}\left(\frac{1+\rho}{1-\rho}\right)
T_n(x),
\label{poschebyexpansion}
\end{equation}
and if $\rho\in(0,1)$ then 
\begin{equation}
\hspace{-0.4cm}\frac{1}
{\rmR^{m}}
T_m
(\zeta_-)
=\frac{1}{(1+\rho)^m}
\sum_{n=0}^\infty\epsilon_n(m)_{2n}{\mathrm P}_{-m}^{-2n}\left(\frac{1-\rho}{1+\rho}\right)
T_n(x).
\label{negchebyexpansion}
\end{equation}
\end{cor}
\noindent {\bf Proof.} 
Using (\ref{gegenposintegral}), 
(\ref{defnChebyshev1intermsofGegenbauer}), and
\begin{equation}
T_m(z)=\sqrt{\frac{\pi}{2}}(z^2-1)^{1/4}P_{m-1/2}^{1/2}(z),
\label{Cheby1defnLegP}
\end{equation}
which follows from 
(\ref{LegendrePdefn}),
(\ref{defnChebyshev1polygauss2F1}), \cite[(15.8.1)]{NIST}.
Analytic continuation to $\rho\in(0,1)$ completes the proof. 
$\hfill\blacksquare$

\medskip

\noindent Note that using 
\cite[(8.6.8)]{Abra}, namely
\[
P_\nu^{1/2}\left(z\right)=
\frac{1}{\sqrt{2\pi}}\left(z^2-1\right)^{-1/4}
\left[\left(z+\sqrt{z^2-1}\right)^{\nu+1/2}+\left(z+\sqrt{z^2-1}\right)^{-\nu-1/2}\right],
\]
and (\ref{Cheby1defnLegP}) one can derive 
the elementary function representation for the 
Chebyshev polynomials of the first kind 
\cite[p.~177]{JakiSharSzab06}.

%
%

\appendix

\section{Definite integrals}
\label{DefiniteintegralsforJacobipolynomials}

As a consequence of the series expansions given above, one may
generate corresponding definite integrals (in a one-step procedure)
as an application of the 
orthogonality relation for these hypergeometric orthogonal polynomials.  
We now describe this correspondence.
Given an expansion over a set of orthogonal polynomials $p_n$ such that
\begin{equation}
f(x)=\sum_{n=0}^\infty a_n p_n(x),
\label{expansion}
\end{equation}
and the orthogonality relation
\begin{equation}
\int_{-1}^{1} p_n(x) p_m(x) w(x) dx = c_n \delta_{n,m},
\label{orthog}
\end{equation}
where $w:(-1,1)\to[0,\infty),$ then using (\ref{orthog}) one has
\[
\int_{-1}^1 f(x) p_n(x) w(x) dx = \sum_{m=0}^\infty a_m \int_{-1}^1 p_n(x) p_m(x) w(x) dx
=a_n c_n,
\]
and therefore
\begin{equation}
a_n=\frac{1}{c_n}\int_{-1}^{1} f(x) p_n(x) w(x) dx.
\label{definitintegral}
\end{equation}
The definite integral expression (\ref{definitintegral}) for the coefficient $a_n$ 
is of equal importance to the 
expansion (\ref{expansion}),
since one may use it to derive the other. 
Integrals of such sort are always of interest since they 
are very likely to find applications in applied mathematics 
and theoretical physics and could be included in tables of integrals
such as \cite{Grad}.

We now give the orthogonality relations for the orthogonal polynomials used in the
main text.  For Jacobi, Gegenbauer, Chebyshev of the second kind, Legendre,
and Chebyshev of the first kind polynomials, the orthogonality relations
can be found in 
\cite[(18.2.1), (18.2.5), Table 18.3.1]{NIST}.
Using the above procedure, we obtain the following definite
integrals for products of special functions with Jacobi, Gegenbauer, Chebyshev
and Legendre polynomials.
Let $m,n\in\N_0$, 
$\alpha,\beta>-1$ such that if $\alpha,\beta\in(-1,0)$ then $\alpha+\beta+1\neq 0$,
$\rho\in\left\{z\in\C:0<|z|<1\right\}\setminus(-1,0]$. Then
\begin{eqnarray*}
&&\hspace{-0.7cm}\int_{-1}^1 \frac{(1-x)^\alpha (1+x)^{\beta/2}}
{\rmR^{\alpha+m+1}}
P_{\alpha+m}^{-\beta}
(\zeta_+)
P_n^{(\alpha,\beta)}(x)dx\\[0.2cm]
&&\hspace{1.1cm}=\frac{2^{\alpha+\beta/2+1}
\Gamma(\alpha+n+1)
(\alpha+\beta+m+1)_{2n}
}
{n!\rho^{(\alpha+1)/2}(1-\rho)^m}P_{-m}^{-\alpha-\beta-2n-1}
\left(\frac{1+\rho}{1-\rho}\right).
\label{defintJac1}
\end{eqnarray*}
Let $\rho\in(0,1).$ Then
\begin{eqnarray*}
&&\hspace{-0.7cm}\int_{-1}^1 \frac{(1-x)^{\alpha/2} (1+x)^{\beta}}
{\rmR^{\beta+m+1}}
{\mathrm P}_{\beta+m}^{-\alpha}
(\zeta_-)
P_n^{(\alpha,\beta)}(x)dx\\[0.2cm]
&&\hspace{1.1cm}=\frac{2^{\alpha/2+\beta+1}
\Gamma(\beta+n+1)
(\alpha+\beta+m+1)_{2n}
}
{n!\rho^{(\beta+1)/2}(1+\rho)^m}{\mathrm P}_{-m}^{-\alpha-\beta-2n-1}
\left(\frac{1-\rho}{1+\rho}\right).
\label{secondintegralcor}
\end{eqnarray*}
Let $\mu\in(-1/2,\infty)\setminus\{0\},$
$\rho\in\{z\in\C:0<|z|<1\}\setminus(-1,0].$ Then
\begin{eqnarray*}
&&\hspace{-0.9cm}\int_{-1}^1 \frac{(1-x^2)^{\mu-1/2}}
{\rmR^{2\mu+m}}
C_m^\mu
(\zeta_+)
C_n^\mu(x)dx
\\[0.2cm]
&&\hspace{2.5cm}=\frac{2^{2-2\mu}\pi\Gamma(2\mu+n)\Gamma(2\mu+2n+m)}
{m!n!\Gamma^2(\mu)\rho^\mu(1-\rho)^m}
P_{-m}^{-2n-2\mu}\left(\frac{1+\rho}{1-\rho}\right).
\end{eqnarray*}
A similar integral on $\rho\in(0,1)$ can be obtained using 
(\ref{negdefintGegenabuer}).

It should be noted that by using
(\ref{UtoC}),
(\ref{PtoC}),
(\ref{defnChebyshev1intermsofGegenbauer}),
the previous definite integral over Gegenbauer polynomials can be written as 
an integral over Chebyshev polynomials of the first and second kind, and well as 
Legendre polynomials.
Let $\lambda\in\C$, $\rho\in\D.$ Then
\begin{eqnarray}
&&\hspace{-0.0cm}\int_{-1}^1 (1-x^2)^{\mu/2-1/4}
P_{\mu-1/2-\lambda}^{1/2-\mu}
\left(\rmR+\rho\right)
{\mathrm P}_{\mu-1/2-\lambda}^{1/2-\mu}
\left(\rmR-\rho\right)
C_n^\mu(x)dx\nonumber\\[0.2cm]
&&\hspace{+5.0cm}=\frac{(\lambda)_n\,(2\mu-\lambda)_n\,\rho^{n+\mu-1/2}\,2^{\mu-1/2}}
{(n+\mu)\Gamma(2\mu)(\mu+1/2)_n\,n!}.\nonumber
\end{eqnarray}
Let $\rho\in(0,1)$. Then 
\begin{eqnarray}
&&\hspace{-0.2cm}\int_{-1}^1\frac{(1-x^2)^{1/4-\mu/2}}
{\rmR^{1/2-\mu+\lambda}}
{\mathrm P}_{\mu-\lambda-1/2}^{1/2-\mu}\left(\frac{1-\rho x}
{\rmR}\right)
C_n^\mu(x)dx\nonumber\\[0.1cm]
&&\hspace{4cm}=\frac{(\lambda)_n\,\sqrt{\pi}\rho^{n+\mu-1/2}\,2^{1/2-\mu}}{(n+\mu)\Gamma(\mu)n!}.\nonumber
\end{eqnarray}
The previous two definite integrals over Gegenbauer polynomials can also be written 
as integrals over Chebyshev polynomials of the second kind and Legendre 
polynomials using (\ref{UtoC}), (\ref{PtoC}).

{\it Acknowledgements}.
The authors thank Hans Volkmer for valuable discussions.
C.~MacKenzie would like to thank the Summer Undergraduate Research Fellowship
program at the National Institute of Standards and Technology for 
financial support while this research was carried out.
We also would like to acknowledge two anonymous referees
whose comments greatly improved the paper.
This work was conducted while H.~S.~Cohl was a National Research Council
Research Postdoctoral Associate in the Applied and Computational
Mathematics Division at the 
National Institute of Standards and Technology, Gaithersburg, Maryland, U.S.A.




\end{document}